\documentclass[12pt]{article}%
\usepackage{amsfonts}
\usepackage{sw20bams}
\usepackage{amsmath}
\usepackage{amssymb}
\usepackage{graphicx}%
\providecommand{\U}[1]{\protect\rule{.1in}{.1in}}
\begin{document}

\title{Variance of Longest Run Duration in a Random Bitstring}
\author{Steven Finch}
\date{July 20, 2020}
\maketitle

\begin{abstract}
We continue an earlier study, starting with unconstrained $n$-bitstrings,
focusing now less on average behavior and more on uncertainty. \ The interplay between

\begin{itemize}
\item longest runs of $0$s and of $1$s, when bitstrings are multus

\item longest runs of $0$s and bitsums (\# of $1$s), when bitstrings are solus

\end{itemize}

\noindent\noindent is examined. \ While negative correlations approach zero as
$n\rightarrow\infty$ in the former (for clumped $1$s), the limit is evidently
nonzero in the latter (for separated $1$s). \ Similar analysis is possible
when both $0$s and $1$s are clumped (bimultus), and when $0$s are clumped but
$1$s are separated (persolus). \ Our methods are experimentally-based.\ 

\end{abstract}

\footnotetext{Copyright \copyright \ 2020 by Steven R. Finch. All rights
reserved.}What can be said about the statistics of the duration $R_{n,1}$ of
the longest run of $1$s in a random bitstring of length $n$? \ Each of the
$2^{n}$ possible bitstrings are assumed to be equally likely; we call this the
\textit{unconstrained} case for reasons that will become clear later. \ The
two summation identities
\[%
%TCIMACRO{\dsum \limits_{j=0}^{\infty}}%
%BeginExpansion
{\displaystyle\sum\limits_{j=0}^{\infty}}
%EndExpansion
j\cdot h_{j}(z)=%
%TCIMACRO{\dsum \limits_{k=0}^{\infty}}%
%BeginExpansion
{\displaystyle\sum\limits_{k=0}^{\infty}}
%EndExpansion
\left(
%TCIMACRO{\dsum \limits_{i=0}^{\infty}}%
%BeginExpansion
{\displaystyle\sum\limits_{i=0}^{\infty}}
%EndExpansion
h_{i}(z)-%
%TCIMACRO{\dsum \limits_{i=0}^{k}}%
%BeginExpansion
{\displaystyle\sum\limits_{i=0}^{k}}
%EndExpansion
h_{i}(z)\right)  ,
\]%
\[%
%TCIMACRO{\dsum \limits_{j=0}^{\infty}}%
%BeginExpansion
{\displaystyle\sum\limits_{j=0}^{\infty}}
%EndExpansion
j^{2}\cdot h_{j}(z)=%
%TCIMACRO{\dsum \limits_{k=0}^{\infty}}%
%BeginExpansion
{\displaystyle\sum\limits_{k=0}^{\infty}}
%EndExpansion
(2k+1)\left(
%TCIMACRO{\dsum \limits_{i=0}^{\infty}}%
%BeginExpansion
{\displaystyle\sum\limits_{i=0}^{\infty}}
%EndExpansion
h_{i}(z)-%
%TCIMACRO{\dsum \limits_{i=0}^{k}}%
%BeginExpansion
{\displaystyle\sum\limits_{i=0}^{k}}
%EndExpansion
h_{i}(z)\right)
\]
give rise to generating functions for the mean and mean square:%
\[
\mathbb{E}(R_{n,1})=\frac{1}{2^{n}}\left[  z^{n}\right]
%TCIMACRO{\dsum \limits_{k=1}^{\infty}}%
%BeginExpansion
{\displaystyle\sum\limits_{k=1}^{\infty}}
%EndExpansion
\left(  \frac{1}{1-2z}-\frac{1-z^{k}}{1-2z+z^{k+1}}\right)  ,
\]%
\[
\mathbb{E}(R_{n,1}^{2})=\frac{1}{2^{n}}\left[  z^{n}\right]
%TCIMACRO{\dsum \limits_{k=1}^{\infty}}%
%BeginExpansion
{\displaystyle\sum\limits_{k=1}^{\infty}}
%EndExpansion
(2k-1)\left(  \frac{1}{1-2z}-\frac{1-z^{k}}{1-2z+z^{k+1}}\right)  .
\]
The expression for $\mathbb{E}(R_{n,1})$ appears in \cite{SF-tcs4}; although a
reference for $\mathbb{E}(R_{n,1}^{2})$ is not known, our expression surely is
not new. \ The Taylor expansion of the numerator series for $\mathbb{E}%
(R_{n,1})$ is \cite{Sn0-tcs4}%
\[
z+4z^{2}+11z^{3}+27z^{4}+62z^{5}+138z^{6}+300z^{7}+643z^{8}+1363z^{9}%
+2866z^{10}+\cdots
\]
and, up to small periodic fluctuations \cite{Byd-tcs4, Sch-tcs4, Alx-tcs4},%
\[
\mathbb{E}(R_{n,1})\sim\frac{\ln(n)}{\ln(2)}-\left(  \frac{3}{2}-\frac{\gamma
}{\ln(2)}\right)
\]
as $n\rightarrow\infty$. \ The Taylor expansion of the numerator series for
$\mathbb{E}(R_{n,1}^{2})$ is%
\[
z+6z^{2}+21z^{3}+61z^{4}+158z^{5}+386z^{6}+902z^{7}+2051z^{8}+4565z^{9}%
+10006z^{10}+
\]
and the variance satisfies%
\[
\mathbb{V}(R_{n,1})\sim\frac{1}{12}+\frac{\pi^{2}}{6\ln(2)^{2}}%
=3.5070480758...
\]
again up to small periodic fluctuations. \ This constant is fairly ubiquitous,
cf. \cite{F1-tcs4, F2-tcs4}. \ Of course, identical results hold for $R_{n,0}%
$, the duration of the longest run of $0$s in a bitstring. \ 

We examine four constrained cases in this paper. \ Let $\Omega$ be a set of
finite bitstrings. \ Rather than exhibit complicated series%
\[
\mathbb{E}(R_{n,1})=\dfrac{1}{d_{n}}\left[  z^{n}\right]  \left\{  G(z)+%
%TCIMACRO{\dsum \limits_{k=1}^{\infty}}%
%BeginExpansion
{\displaystyle\sum\limits_{k=1}^{\infty}}
%EndExpansion
\left(  H(z)-H_{k}(z)\right)  \right\}  ,
\]%
\[
\mathbb{E}(R_{n,1}^{2})=\dfrac{1}{d_{n}}\left[  z^{n}\right]  \left\{  G(z)+%
%TCIMACRO{\dsum \limits_{k=1}^{\infty}}%
%BeginExpansion
{\displaystyle\sum\limits_{k=1}^{\infty}}
%EndExpansion
(2k-1)\left(  H(z)-H_{k}(z)\right)  \right\}  ,
\]
for each choice of $\Omega$, we instead simply provide the denominator%
\[
d_{n}=\text{the count of all }\omega\in\Omega\text{ of length }n
\]
as well as the generating function $H_{k}(z)$ for the count of all $\omega
\in\Omega$ with no runs of $k$ $1$s (or $k$ $0$s, depending on the scenario).
A nonzero correction term $G(z)$ is often required too. For the preceding,
clearly $d_{n}=2^{n}$, $G(z)=0$ and
\[%
\begin{array}
[c]{ccc}%
H_{k}(z)=\dfrac{1-z^{k}}{1-2z+z^{k+1}}, &  & H(z)=\lim\limits_{k\rightarrow
\infty}H_{k}(z)=\dfrac{1}{1-2z}.
\end{array}
\]
For $\mathbb{E}(R_{n,1}^{3})$ and $\mathbb{E}(R_{n,1}^{4})$, the factor $2k-1$
in the $\mathbb{E}(R_{n,1}^{2})$ formula would be replaced by $3k^{2}-3k+1$
and $4k^{3}-6k^{2}+4k-1$ respectively.

\section{Constrained Cases}

The first two examples of constrained bitstrings were introduced in
\cite{Fi-tcs4}.

\begin{itemize}
\item A bitstring is \textbf{solus} if all of its $1$s are isolated.

\item A bitstring is \textbf{multus} if each of its $1$s possess at least one
neighboring $1$.
\end{itemize}

\noindent Counts of solus $n$-bitstrings have a quadratic character:%
\[%
%TCIMACRO{\dsum \limits_{n=0}^{\infty}}%
%BeginExpansion
{\displaystyle\sum\limits_{n=0}^{\infty}}
%EndExpansion
d_{n}z^{n}=\frac{1+z}{1-z-z^{2}}=1+2z+3z^{2}+5z^{3}+8z^{4}+13z^{5}%
+21z^{6}+34z^{7}+\cdots,
\]
whereas counts of multus $n$-bitstrings have a cubic character:%
\[%
%TCIMACRO{\dsum \limits_{n=0}^{\infty}}%
%BeginExpansion
{\displaystyle\sum\limits_{n=0}^{\infty}}
%EndExpansion
d_{n}z^{n}=\frac{1-z+z^{2}}{1-2z+z^{2}-z^{3}}=1+z+2z^{2}+4z^{3}+7z^{4}%
+12z^{5}+21z^{6}+37z^{7}+\cdots.
\]

\noindent The remaining two examples are new, as far as is known.

\begin{itemize}
\item A bitstring is \textbf{bimultus} if each of its $1$s possess at least
one neighboring $1$ and each of its $0$s possess at least one neighboring $0$.
\ Both isolated $1$ bits and isolated $0$ bits are avoided in such bitstrings;
a certain $0\leftrightarrow1$ symmetry holds here.

\item A bitstring is \textbf{persolus} if all of its $1$s are isolated and
each of its $0$s possess at least one neighboring $0$. \ That is, while $1$s
in solus bitstrings are alone, $1$s in persolus bitstrings are \textit{very} alone.
\end{itemize}

\noindent Counts of bimultus $n$-bitstrings have a quadratic character
\cite{Sn1-tcs4}%
\[%
%TCIMACRO{\dsum \limits_{n=0}^{\infty}}%
%BeginExpansion
{\displaystyle\sum\limits_{n=0}^{\infty}}
%EndExpansion
d_{n}z^{n}=\frac{2z^{2}}{1-z-z^{2}}=2z^{2}+2z^{3}+4z^{4}+6z^{5}+10z^{6}%
+16z^{7}+\cdots,
\]
whereas counts of persolus $n$-bitstrings have a cubic character%
\[%
%TCIMACRO{\dsum \limits_{n=0}^{\infty}}%
%BeginExpansion
{\displaystyle\sum\limits_{n=0}^{\infty}}
%EndExpansion
d_{n}z^{n}=\frac{z\left(  1+2z^{2}\right)  }{1-z-z^{3}}=z+z^{2}+3z^{3}%
+4z^{4}+5z^{5}+8z^{6}+12z^{7}+\cdots.
\]

\section{Bitsums}

Given a set $\Omega$ of finite bitstrings, what can be said about the bitsum
$S_{n}$ of a random $\omega\in\Omega$ of length $n$? \ If $\Omega$ is
unconstrained, i.e., if all $2^{n}$ strings are included in the sample, then%
\[%
\begin{array}
[c]{ccc}%
\mathbb{E}(S_{n})=n/2, &  & \mathbb{V}(S_{n})=n/4
\end{array}
\]
because a sum of $n$ independent Bernoulli($1/2$) variables is Binomial($n$%
,$1/2$). \ Expressed differently, the average density of $1$s in a random
unconstrained string is $1/2$, with a corresponding variance $1/4$. \ 

We previously covered solus and multus bitstrings in \cite{Fi-tcs4}. \ If
$\Omega$ consists of bimultus bitstrings, then the total bitsum $a_{n}$ of all
$\omega\in\Omega$ of length $n$ has generating function \cite{Sn1-tcs4}
\[%
%TCIMACRO{\dsum \limits_{n=0}^{\infty}}%
%BeginExpansion
{\displaystyle\sum\limits_{n=0}^{\infty}}
%EndExpansion
a_{n}z^{n}=\frac{z^{2}(2-z)}{\left(  1-z-z^{2}\right)  ^{2}}=2z^{2}%
+3z^{3}+8z^{4}+15z^{5}+30z^{6}+\cdots
\]
and the total bitsum squared $b_{n}$ has generating function%
\[%
%TCIMACRO{\dsum \limits_{n=0}^{\infty}}%
%BeginExpansion
{\displaystyle\sum\limits_{n=0}^{\infty}}
%EndExpansion
b_{n}z^{n}=\frac{z^{2}\left(  4-7z+4z^{2}-z^{3}+4z^{4}-z^{5}\right)  }{\left(
1-z+z^{2}\right)  \left(  1-z-z^{2}\right)  ^{3}}=4z^{2}+9z^{3}+24z^{4}%
+51z^{5}+114z^{6}+\cdots;
\]
hence $c_{n}=d_{n}b_{n}-a_{n}^{2}$ has generating function
\begin{align*}%
%TCIMACRO{\dsum \limits_{n=0}^{\infty}}%
%BeginExpansion
{\displaystyle\sum\limits_{n=0}^{\infty}}
%EndExpansion
c_{n}z^{n}  &  =\frac{z^{2}\left(  4-11z+11z^{2}-13z^{3}+2z^{4}+17z^{5}%
-5z^{6}-z^{7}\right)  }{(1+z)^{2}\left(  1-3z+z^{2}\right)  ^{2}\left(
1-z+2z^{2}+z^{3}+z^{4}\right)  }\\
&  =4z^{2}+9z^{3}+32z^{4}+81z^{5}+240z^{6}+\cdots.
\end{align*}
Standard techniques \cite{SF-tcs4} give asymptotics \
\[
\lim_{n\rightarrow\infty}\frac{\mathbb{E}(S_{n})}{n}=\lim_{n\rightarrow\infty
}\frac{a_{n}}{nd_{n}}=\frac{1}{2},
\]%
\[
\lim_{n\rightarrow\infty}\frac{\mathbb{V}(S_{n})}{n}=\lim_{n\rightarrow\infty
}\frac{c_{n}}{nd_{n}^{2}}=\frac{5+3\sqrt{5}}{40}=0.2927050983...
\]
for the average density of $1$s in a random bimultus string and corresponding
variance. \ 

If $\Omega$ instead consists of persolus bitstrings, then the total bitsum
$a_{n}$ of all $\omega\in\Omega$ of length $n$ has generating function
\cite{Sn1-tcs4}%
\[%
%TCIMACRO{\dsum \limits_{n=0}^{\infty}}%
%BeginExpansion
{\displaystyle\sum\limits_{n=0}^{\infty}}
%EndExpansion
a_{n}z^{n}=\frac{z(1-z+z^{2})^{2}}{\left(  1-z-z^{3}\right)  ^{2}}%
=z+2z^{3}+4z^{4}+5z^{5}+10z^{6}+\cdots
\]
and the total bitsum squared $b_{n}$ has generating function%
\[%
%TCIMACRO{\dsum \limits_{n=0}^{\infty}}%
%BeginExpansion
{\displaystyle\sum\limits_{n=0}^{\infty}}
%EndExpansion
b_{n}z^{n}=\frac{z\left(  1-z+z^{2}\right)  ^{2}\left(  1-z+z^{3}\right)
}{\left(  1-z-z^{3}\right)  ^{3}}=z+2z^{3}+6z^{4}+7z^{5}+16z^{6}+\cdots;
\]
hence $c_{n}=d_{n}b_{n}-a_{n}^{2}$ has generating function
\begin{align*}%
%TCIMACRO{\dsum \limits_{n=0}^{\infty}}%
%BeginExpansion
{\displaystyle\sum\limits_{n=0}^{\infty}}
%EndExpansion
c_{n}z^{n}  &  =\frac{z^{3}\left(  2+4z-6z^{2}-6z^{3}-16z^{4}-8z^{5}%
+8z^{6}+14z^{7}+5z^{8}-2z^{9}-3z^{10}-z^{11}\right)  }{\left(  1-z-2z^{2}%
-z^{3}\right)  ^{2}\left(  1+z^{2}-z^{3}\right)  ^{3}}\\
&  =2z^{3}+8z^{4}+10z^{5}+28z^{6}+\cdots.
\end{align*}
We obtain asymptotics \
\begin{align*}
\lim_{n\rightarrow\infty}\frac{\mathbb{E}(S_{n})}{n}  &  =\lim_{n\rightarrow
\infty}\frac{a_{n}}{nd_{n}}\\
&  =\frac{1}{3}\left[  1-\left(  \frac{31+3\sqrt{93}}{1922}\right)
^{1/3}-\left(  \frac{31-3\sqrt{93}}{1922}\right)  ^{1/3}\right] \\
&  =0.1942540040...,
\end{align*}%
\begin{align*}
\lim_{n\rightarrow\infty}\frac{\mathbb{V}(S_{n})}{n}  &  =\lim_{n\rightarrow
\infty}\frac{c_{n}}{nd_{n}^{2}}\\
&  =\frac{1}{2883}\left(  \frac{93}{2}\right)  ^{1/3}\left[  \left(
8649+457\sqrt{93}\right)  ^{1/3}+\left(  8649-457\sqrt{93}\right)
^{1/3}\right] \\
&  =0.0495615175...
\end{align*}
for the average density of $1$s in a random persolus string and corresponding
variance. \ Unsurprisingly $0.588>1/2>0.276>0.194$ and
$0.292>0.281>1/4>0.089>0.049$, where $\{0.276,0.089\}$ are estimates
associated with solus strings and $\{0.588,0.281\}$ are estimates associated
with multus strings \cite{Fi-tcs4}. \ While insisting on $0\leftrightarrow1$
symmetry forces equiprobability, it also increases the variance, but only slightly.

\section{Longest Bitruns}

Given a set $\Omega$ of finite bitstrings, what can be said about the duration
$R_{n,1}$ of the longest run of $1$s in a random $\omega\in\Omega$ of length
$n$? \ We have already discussed the case when $\Omega$ is unconstrained.
\ Preliminary coverage for constrained $\Omega$ (for means, but not mean
squares) occurred in \cite{Fi-tcs4}.

If $\Omega$ consists of solus bitstrings, then it makes little sense to talk
about $1$-runs. \ For $0$-runs, over all $\omega\in\Omega$, we have%
\[%
\begin{array}
[c]{ccc}%
H_{k}(z)=\dfrac{1+z-z^{k}-z^{k+1}}{1-z-z^{2}+z^{k+1}}, &  & H(z)=\dfrac
{1+z}{1-z-z^{2}};
\end{array}
\]
the Taylor expansion of the numerator series for $\mathbb{E}(R_{n,0})$ is
\[
z+4z^{2}+9z^{3}+18z^{4}+34z^{5}+62z^{6}+110z^{7}+192z^{8}+331z^{9}%
+565z^{10}+\cdots
\]
and the Taylor expansion of the numerator series for $\mathbb{E}(R_{n,0}^{2})$
is
\[
z+6z^{2}+19z^{3}+48z^{4}+106z^{5}+218z^{6}+424z^{7}+798z^{8}+1463z^{9}%
+2631z^{10}+\cdots.
\]
Let us abbreviate such series as $\operatorname{num}_{0}^{1}$ and
$\operatorname{num}_{0}^{2}$ for simplicity -- likewise $\operatorname{num}%
_{1}^{1}$ and $\operatorname{num}_{1}^{2}$ -- and let $\varphi=(1+\sqrt
{5})/2=1.6180339887...$ denote the Golden mean. \ It is conjectured that, up
to small periodic fluctuations,%
\[
\mathbb{E}(R_{n,0})\sim\frac{\ln(n)}{\ln(\varphi)}-\left(  2-\frac{\gamma}%
{\ln(\varphi)}\right)  ,
\]%
\[
\mathbb{V}(R_{n,0})\sim\frac{1}{12}+\frac{\pi^{2}}{6\ln(\varphi)^{2}%
}=7.1868910445...
\]
as $n\rightarrow\infty$. \ 

If $\Omega$ consists of multus bitstrings, then we can talk both about
$1$-runs:%
\[
G(z)=\frac{-z}{\left(  1-z\right)  \left(  1-z+z^{2}\right)  },
\]%
\[%
\begin{array}
[c]{ccc}%
H_{k}(z)=\dfrac{1+z^{2}-z^{k-1}-z^{k}}{1-2z+z^{2}-z^{3}+z^{k+1}}z, &  &
H(z)=\dfrac{1+z^{2}}{1-2z+z^{2}-z^{3}}z;
\end{array}
\]%
\[
\operatorname{num}_{1}^{1}=2z^{2}+7z^{3}+16z^{4}+32z^{5}+62z^{6}%
+118z^{7}+221z^{8}+409z^{9}+751z^{10}+\cdots
\]%
\[
\operatorname{num}_{1}^{2}=4z^{2}+17z^{3}+46z^{4}+104z^{5}+220z^{6}%
+448z^{7}+889z^{8}+1729z^{9}+3313z^{10}+\cdots
\]
and $0$-runs:%
\[
G(z)=0,
\]%
\[%
\begin{array}
[c]{ccc}%
H_{k}(z)=\dfrac{1+z^{2}-z^{k-1}+z^{k}-2z^{k+1}}{1-2z+z^{2}-z^{3}+z^{k+2}}z, &
& H(z)=\dfrac{1+z^{2}}{1-2z+z^{2}-z^{3}}z;
\end{array}
\]%
\[
\operatorname{num}_{0}^{1}=z+2z^{2}+5z^{3}+11z^{4}+23z^{5}+45z^{6}%
+87z^{7}+165z^{8}+309z^{9}+573z^{10}+\cdots,
\]%
\[
\operatorname{num}_{0}^{2}=z+4z^{2}+11z^{3}+27z^{4}+63z^{5}+135z^{6}%
+281z^{7}+565z^{8}+1115z^{9}+2161z^{10}+\cdots.
\]
Letting%
\[
\psi=\frac{1}{3}\left[  2+\left(  \frac{25+3\sqrt{69}}{2}\right)
^{1/3}+\left(  \frac{25-3\sqrt{69}}{2}\right)  ^{1/3}\right]
=1.7548776662...,
\]
it is conjectured that%
\[
\mathbb{E}(R_{n,1})\sim\mathbb{E}(R_{n,0})+1\sim\frac{\ln(n)}{\ln(\psi
)}-\left(  \frac{3}{2}-\frac{\gamma}{\ln(\psi)}\right)  ,
\]%
\[
\mathbb{V}(R_{n,1})\sim\mathbb{V}(R_{n,0})\sim\frac{1}{12}+\frac{\pi^{2}}%
{6\ln(\psi)^{2}}=5.2840019997...
\]
as $n\rightarrow\infty$. \ 

If $\Omega$ consists of bimultus strings, there is symmetry (just as for the
unconstrained case). \ We have
\[
G(z)=\frac{-z(1-z+z^{2})^{2}}{(1-z)(1-z+z^{3})},
\]%
\[%
\begin{array}
[c]{ccc}%
H_{k}(z)=\dfrac{2-2z+2z^{2}-z^{k-2}+z^{k-1}-2z^{k}}{1-2z+z^{2}-z^{4}+z^{k+2}%
}z^{2}, &  & H(z)=\dfrac{2-2z+2z^{2}}{1-2z+z^{2}-z^{4}}z^{2};
\end{array}
\]%
\[
\operatorname{num}_{0}^{1}=\operatorname{num}_{1}^{1}=2z^{2}+3z^{3}%
+8z^{4}+15z^{5}+28z^{6}+50z^{7}+87z^{8}+150z^{9}+255z^{10}+\cdots,
\]%
\[
\operatorname{num}_{0}^{2}=\operatorname{num}_{1}^{2}=4z^{2}+9z^{3}%
+24z^{4}+51z^{5}+102z^{6}+196z^{7}+361z^{8}+656z^{9}+1165z^{10}+\cdots.
\]
It is conjectured that, up to small periodic fluctuations,%
\[
\mathbb{E}(R_{n,0})\sim\mathbb{E}(R_{n,1})\sim\frac{\ln(n)}{\ln(\varphi
)}-\left(  \frac{5}{2}-\frac{\gamma}{\ln(\varphi)}\right)  ,
\]%
\[
\mathbb{V}(R_{n,0})\sim\mathbb{V}(R_{n,1})\sim\frac{1}{12}+\frac{\pi^{2}}%
{6\ln(\varphi)^{2}}=7.1868910445...
\]
as $n\rightarrow\infty$. \ The same asymptotic variance occurred for solus bitstrings.

If $\Omega$ consists of persolus strings, then it makes little sense to talk
about $1$-runs. \ For $0$-runs, we have%

\[
G(z)=\frac{-z(1+z)^{2}}{1+z^{2}},
\]%
\[%
\begin{array}
[c]{ccc}%
H_{k}(z)=\dfrac{1+2z^{2}-z^{k-1}-2z^{k}}{1-z-z^{3}+z^{k+1}}z, &  &
H(z)=\dfrac{1+2z^{2}}{1-z-z^{3}}z;
\end{array}
\]%
\[
\operatorname{num}_{0}^{1}=2z^{2}+7z^{3}+12z^{4}+18z^{5}+30z^{6}%
+49z^{7}+76z^{8}+118z^{9}+183z^{10}+\cdots,
\]%
\[
\operatorname{num}_{0}^{2}=4z^{2}+17z^{3}+38z^{4}+70z^{5}+128z^{6}%
+227z^{7}+384z^{8}+636z^{9}+1037z^{10}+\cdots.
\]
Letting%
\[
\theta=\frac{1}{3}\left[  1+\left(  \frac{29+3\sqrt{93}}{2}\right)
^{1/3}+\left(  \frac{29-3\sqrt{93}}{2}\right)  ^{1/3}\right]
=1.4655712318...,
\]
it is conjectured that%
\[
\mathbb{E}(R_{n,0})\sim\frac{\ln(n)}{\ln(\theta)}-\left(  \frac{5}{2}%
-\frac{\gamma}{\ln(\theta)}\right)  ,
\]%
\[
\mathbb{V}(R_{n,0})\sim\frac{1}{12}+\frac{\pi^{2}}{6\ln(\theta)^{2}%
}=11.3414222234...
\]
as $n\rightarrow\infty$. \ The constants $\varphi$ and $\psi$ also appeared in
\cite{Fi-tcs4}; $\theta$ was called Moore's constant in \cite{Fn-tcs4} and is
the limit of a certain fundamental iteration.

\section{Cross-Covariances I}

Let us return to the unconstrained case. \ Exhibiting $\mathbb{E}%
(R_{n,0}R_{n,1})$ in a manner parallel to old formulas in our introduction
seems impossible: no analogous summation identity for $%
%TCIMACRO{\tsum \nolimits_{i=0}^{\infty}}%
%BeginExpansion
{\textstyle\sum\nolimits_{i=0}^{\infty}}
%EndExpansion%
%TCIMACRO{\tsum \nolimits_{j=0}^{\infty}}%
%BeginExpansion
{\textstyle\sum\nolimits_{j=0}^{\infty}}
%EndExpansion
i\,j\cdot h_{i,j}(z)$ apparently exists. \ Thus new formulas are somewhat less
tidy, but nevertheless workable. \ The number of bitstrings with no runs of
$i$ $1$s and no runs of $j$ $0$s has generating function%
\[
f_{i,j}(z)=\frac{1-z^{i}-z^{j}+z^{i+j}}{1-2z+z^{i+1}+z^{j+1}-z^{i+j}}%
\]
hence%
\[
\mathbb{E}(R_{n,0}R_{n,1})=\frac{1}{d_{n}}\left[  z^{n}\right]
%TCIMACRO{\dsum \limits_{i=1}^{\infty}}%
%BeginExpansion
{\displaystyle\sum\limits_{i=1}^{\infty}}
%EndExpansion%
%TCIMACRO{\dsum \limits_{j=1}^{\infty}}%
%BeginExpansion
{\displaystyle\sum\limits_{j=1}^{\infty}}
%EndExpansion
i\,j\left\{  f_{i+1,j+1}(z)-f_{i,j+1}(z)-f_{i+1,j}(z)+f_{i,j}(z)\right\}
\]
where $d_{n}=2^{n}$. \ The Taylor expansion of the numerator series for
$\mathbb{E}(R_{n,0}R_{n,1})$ is
\[
2z^{2}+10z^{3}+34z^{4}+96z^{5}+248z^{6}+604z^{7}+1418z^{8}+3240z^{9}%
+7260z^{10}+\cdots
\]
and the correlation coefficient
\[
\rho=\frac{\mathbb{E}(R_{n,0}R_{n,1})-\mathbb{E}(R_{n,0})\mathbb{E}(R_{n,1}%
)}{\sqrt{\mathbb{V}(R_{n,0})}\sqrt{\mathbb{V}(R_{n,1})}}%
\]
is prescribed numerically in Table 1 for $n=10,20,\ldots,70$. \ These results
complement those in \cite{PM-tcs4}. \ 

For multus bitstrings, since $1$s are clumped (but $0$s are not necessarily
so), the associated generating function
\[
f_{i,j}(z)=\frac{1+z^{2}-z^{i-1}-z^{i}-z^{j-1}+z^{j}-2z^{j+1}+2z^{i+j-1}%
}{1-2z+z^{2}-z^{3}+z^{i+1}+z^{j+2}-z^{i+j}}z
\]
is unsurprisingly asymmetric in $i$ and $j$. \ The associated Taylor expansion
is%
\[
4z^{3}+16z^{4}+45z^{5}+106z^{6}+232z^{7}+484z^{8}+977z^{9}+1927z^{10}+\cdots
\]
and, again, the corresponding $\rho$ is prescribed in Table 1. \ Correlations
are all negative but approach zero as $n$ increases. \ We observe a slightly
stronger dependency between $R_{n,0}$ and $R_{n,1}$ for multus strings than
for unconstrained strings. Calculating $f_{i,j}(z)$ for bimultus strings
remains open. \ Simulation suggests that dependence is greater still for the
bimultus case.

\begin{center}%
\[%
\begin{tabular}
[c]{|c|c|c|}\hline
$n$ & $\rho$ for unconstrained case & $\rho$ for multus case\\\hline
$10$ & $-0.383683$ & $-0.443900$\\\hline
$20$ & $-0.225906$ & $-0.256080$\\\hline
$30$ & $-0.165175$ & $-0.187941$\\\hline
$40$ & $-0.132345$ & $-0.151033$\\\hline
$50$ & $-0.111286$ & $-0.127411$\\\hline
$60$ & $-0.096550$ & $-0.110810$\\\hline
$70$ & $-0.085616$ & $-0.098434$\\\hline
\end{tabular}
\ \
\]
Table 1:\ Correlation $\rho(R_{n,0},R_{n,1})$ as a function of $n$.\bigskip
\end{center}

\section{Cross-Covariances II}

Let us return to the solus case. \ Being isolated, each $1$ acts as barrier to
gathering $0$s; we wonder to what extent the (random) number of such walls
affects the largest crowd size. \ To calculate $\mathbb{E}(R_{n,0}S_{n})$
seems to be difficult. \ The number of bitstrings with less than two $1$s and
no runs of $k$ $0$s has generating function%
\[
f_{2,k}(z)=%
%TCIMACRO{\dsum \limits_{n=1}^{\infty}}%
%BeginExpansion
{\displaystyle\sum\limits_{n=1}^{\infty}}
%EndExpansion
a_{n}z^{n}%
\]
(a polynomial!)\ with%
\[
a_{n}=\left\{
\begin{array}
[c]{lll}%
n+1 &  & \text{if }1\leq n\leq k-1\text{,}\\
2k-n &  & \text{if }k\leq n\leq2k-1\text{,}\\
0 &  & \text{otherwise}%
\end{array}
\right.
\]
assuming $k\geq2$. \ For example,%
\[
\{a_{n}\}_{n=1}^{2k-1}=\left\{  2,3,4,5,6,7,7,6,5,4,3,2,1\right\}
\]
when $k=7$. \ 

The number of bitstrings with less than three $1$s and no runs of $k$ $0$s has
generating function $f_{3,k}(z)$ with%
\[
a_{n}=\left\{
\begin{array}
[c]{lll}%
2-\dfrac{n}{2}+\dfrac{n^{2}}{2} &  & \text{if }1\leq n\leq k-1\text{,}\\
a_{n-1}+k-2 &  & \text{if }k\leq n\leq k+2\text{,}\\
a_{n-1}+3k-2n+2 &  & \text{if }k+3\leq n\leq2k\text{,}\\
\dfrac{3k-n}{2}+\dfrac{(3k-n)^{2}}{2} &  & \text{if }2k+1\leq n\leq
3k-1\text{,}\\
0 &  & \text{otherwise}%
\end{array}
\right.
\]
assuming $k\geq2$. \ For example,%
\[
\{a_{n}\}_{n=1}^{3k-1}=\left\{
2,3,5,8,12,17,22,27,32,35,36,35,32,27,21,15,10,6,3,1\right\}
\]
when $k=7$. \ 

The number of bitstrings with less than four $1$s and no runs of $k$ $0$s has
generating function $f_{4,k}(z)$ with%
\[
a_{n}=\left\{
\begin{array}
[c]{lll}%
2 &  & \text{if }n=1\text{,}\\
1-\delta_{k,n}+\dfrac{7(n-1)}{3}-\dfrac{(n-1)^{2}}{2}+\dfrac{(n-1)^{3}}{6} &
& \text{if }2\leq n\leq k\text{,}\\
a_{n-1}+u(k,n) &  & \text{if }k+1\leq n\leq2k+2\text{,}\\
a_{n-1}-v(k,n) &  & \text{if }2k+3\leq n\leq3k\text{,}\\
\dfrac{4k-n}{3}+\dfrac{(4k-n)^{2}}{2}+\dfrac{(4k-n)^{3}}{6} &  & \text{if
}3k+1\leq n\leq4k-1\text{,}\\
0 &  & \text{otherwise}%
\end{array}
\right.
\]
assuming $k\geq2$, where%
\[
u(k,n)=\left\{
\begin{array}
[c]{lll}%
\dfrac{-k^{2}+(2n-5)k-w(n-k)}{2} &  & \text{if }k+1\leq n\leq2k\text{,}\\
2(n-k-3) &  & \text{if }n=2k+1\text{,}\\
0 &  & \text{otherwise;}%
\end{array}
\right.
\]%
\[
v(k,n)=\dfrac{-20k^{2}+(16n-30)k-\left(  3n^{2}-11n+12\right)  }{2}\text{;}%
\]%
\[
w(m)=\left\{
\begin{array}
[c]{lll}%
-2 &  & \text{if }m=1\text{,}\\
2 &  & \text{if }m=2\text{,}\\
3m^{2}-13m+20 &  & \text{if }m\geq3\text{;}%
\end{array}
\right.
\]
$\delta_{k,n}$ is $1$ when $k=n$ and is $0$ otherwise. \ For example,%
\begin{align*}
\{a_{n}\}_{n=1}^{4k-1}  &  =\left\{
2,3,5,8,13,21,32,47,67,91,118,145,169,187,\right. \\
&  \left.  \;\;\;\;\;197,197,186,166,140,111,82,56,35,20,10,4,1\right\}
\end{align*}
when $k=7$. \ 

An expression for $f_{5,k}(z)$, the generating function corresponding to
bitstrings with less than five $1$s and no runs of $k$ $0$s, exists but awaits
simplication. \ For example,%
\begin{align*}
\{a_{n}\}_{n=1}^{5k-1}  &  =\left\{
2,3,5,8,13,21,33,52,82,126,188,271,376,500,637,777,907,1013,\right. \\
&  \left.
\;\;\;\;\;1081,1102,1073,997,882,741,590,444,314,207,126,70,35,15,5,1\right\}
\end{align*}
when $k=7$. \ For arbitrary $k\geq2$, clearly $a_{n}$ is equal to%
\[
4-\delta_{k,n}-\dfrac{n-2}{4}+\dfrac{35(n-2)^{2}}{24}-\dfrac{(n-2)^{3}}%
{4}+\dfrac{(n-2)^{4}}{24}%
\]
for $3\leq n\leq k$ and%
\[
\dfrac{5k-n}{4}+\dfrac{11(5k-n)^{2}}{24}+\dfrac{(5k-n)^{3}}{4}+\dfrac
{(5k-n)^{4}}{24}%
\]
for $4k+1\leq n\leq5k-1$. \ Also,\bigskip\
\[
a_{n}=\left\{
\begin{array}
[c]{lll}%
a_{n-1}+u(k,n) &  & \text{if }k+1\leq n\leq2k+1\text{,}\\
a_{n-1}-v(k,n) &  & \text{if }3k+2\leq n\leq4k
\end{array}
\right.
\]
where\bigskip%
\[
u(k,n)=\delta_{2k+1,n}+\dfrac{k^{3}-(3n-12)k^{2}+(3n^{2}-24n+59)k-w(n-k)}{6};
\]
\bigskip%
\[
v(k,n)=3\delta_{3k+2,n}+\dfrac{-195k^{3}+(165n-426)k^{2}-\left(
45n^{2}-228n+309\right)  k+\left(  4n^{3}-30n^{2}+80n-72\right)  }{6}\text{;}%
\]
\bigskip%
\[
w(m)=\left\{
\begin{array}
[c]{lll}%
54 &  & \text{if }m=1\text{,}\\
30 &  & \text{if }2\leq m\leq3\text{,}\\
4m^{3}-42m^{2}+176m-240 &  & \text{if }m\geq4\text{;}%
\end{array}
\right.
\]
\bigskip%
\[
a_{3k+1}=\dfrac{11k^{4}-2k^{3}-35k^{2}-22k+72}{24};
\]
\bigskip%
\[
\left(  \text{index of }\max_{1\leq n\leq5k-1}a_{n}\right)  =\left\{
\begin{array}
[c]{lll}%
(5k+5)/2 &  & \text{if }k\geq3\text{ is odd,}\\
(5k+4)/2 &  & \text{if }k\geq4\text{ is even;}%
\end{array}
\right.
\]%
\[
\left(  \max_{1\leq n\leq5k-1}a_{n}\right)  =\left\{
\begin{array}
[c]{lll}%
\left(  115k^{4}-184k^{3}-22k^{2}-104k+387\right)  /192 &  & \text{if }%
k\geq3\text{ is odd,}\\
\left(  115k^{4}-184k^{3}-52k^{2}+16k+192\right)  /192 &  & \text{if }%
k\geq2\text{ is even.}%
\end{array}
\right.  \bigskip
\]
\noindent The interval $2k+2\leq n\leq3k$ deserves more attention. \ From the
approximation
\begin{align*}
&
%TCIMACRO{\dsum \limits_{k=2}^{\infty}}%
%BeginExpansion
{\displaystyle\sum\limits_{k=2}^{\infty}}
%EndExpansion
k\left\{  f_{2,k+1}(z)-f_{2,k}(z)\right\}  +2%
%TCIMACRO{\dsum \limits_{k=2}^{\infty}}%
%BeginExpansion
{\displaystyle\sum\limits_{k=2}^{\infty}}
%EndExpansion
k\left\{  f_{3,k+1}(z)-f_{2,k+1}(z)-f_{3,k}(z)+f_{2,k}(z)\right\} \\
&  +3%
%TCIMACRO{\dsum \limits_{k=2}^{\infty}}%
%BeginExpansion
{\displaystyle\sum\limits_{k=2}^{\infty}}
%EndExpansion
k\left\{  f_{4,k+1}(z)-f_{3,k+1}(z)-f_{4,k}(z)+f_{3,k}(z)\right\}  +4%
%TCIMACRO{\dsum \limits_{k=2}^{\infty}}%
%BeginExpansion
{\displaystyle\sum\limits_{k=2}^{\infty}}
%EndExpansion
k\left\{  f_{5,k+1}(z)-f_{4,k+1}(z)-f_{5,k}(z)+f_{4,k}(z)\right\}
\end{align*}
we obtain the Taylor expansion of the numerator series for $\mathbb{E}%
(R_{n,0}S_{n})$:%
\[
2z^{2}+7z^{3}+18z^{4}+43z^{5}+94z^{6}+196z^{7}+392z^{8}+764z^{9}%
+1454z^{10}+\cdots
\]
(every exhibited coefficient is correct). \ We would need $f_{6,k}(z)$,
$f_{7,k}(z)$, \ldots\ to achieve the precision necessary to adequately
estimate $\rho(R_{n,0},S_{n})$ for large $n$. \ Simulation suggests that
correlations are all negative but, unlike the previous section, tend to a
nonzero quantity (possibly $<-1/10$?) as $n$ approaches infinity. \ We have
not yet attempted to study the persolus case.

\section{Cross-Covariances III}

This section is an addendum to the preceding. \ A recent paper \cite{NS-tcs4}
gave an impressive recursion for the number $F_{n}(x,y)$ of unconstrained
bitstrings of length $n$ containing exactly $x$ $0$s and a longest run of
exactly $y$ $0$s:
\[
F_{n}(x,y)=\left\{
\begin{array}
[c]{lll}%
%TCIMACRO{\dsum \limits_{i=\kappa}^{y-1}}%
%BeginExpansion
{\displaystyle\sum\limits_{i=\kappa}^{y-1}}
%EndExpansion
F_{n-i-1}(x-i,y)+%
%TCIMACRO{\dsum \limits_{j=0}^{y}}%
%BeginExpansion
{\displaystyle\sum\limits_{j=0}^{y}}
%EndExpansion
F_{n-y-1}(x-y,j) &  & \text{if }1\leq x\leq n-2\text{ and }\varepsilon
_{n}(x,y)=1,\\
\lambda_{n}(y) &  & \text{if }x=n-1\text{ and }\varepsilon_{n}(x,y)=1,\\
0 &  & \text{otherwise,}%
\end{array}
\right.
\]%
\[%
\begin{array}
[c]{ccc}%
F_{n}(0,0)=1-\kappa, &  & F_{n}(n,n)=1
\end{array}
\]
where $\kappa=0$,%
\[
\varepsilon_{n}(x,y)=\left\{
\begin{array}
[c]{lll}%
1 &  & \text{if }n\geq2\text{ and }\left\lfloor \dfrac{n}{n-x+1}\right\rfloor
\leq y\leq x,\\
0 &  & \text{otherwise;}%
\end{array}
\right.
\]%
\[
\lambda_{n}(y)=\left\{
\begin{array}
[c]{lll}%
1 &  & \text{if }n\text{ is odd and }y=\dfrac{n-1}{2},\\
2 &  & \text{otherwise.}%
\end{array}
\right.
\]
By a similar argument, we deduce the number $\widetilde{F}_{n}(x,y)$ of solus
bitstrings of length $n$ containing exactly $x$ $0$s and a longest run of
exactly $y$ $0$s:%
\[
\widetilde{F}_{n}(x,y)=\left\{
\begin{array}
[c]{lll}%
F_{n-1}(x,y)+F_{n}(x,y) &  & \text{if }n\geq2,\\
\delta_{x,y} &  & \text{otherwise}%
\end{array}
\right.
\]
where $F_{n}(x,y)$ is defined recursively as before, with the same
$\varepsilon_{n}(x,y)$ but with $\kappa=1$ and a different $\lambda_{n}(y)$:
\[
\lambda_{n}(y)=\left\{
\begin{array}
[c]{lll}%
\left[
\begin{array}
[c]{lll}%
1 &  & \text{if }y=\dfrac{n-1}{2}\text{ or }y=n-1,\\
2 &  & \text{otherwise}%
\end{array}
\right.  &  & \text{if }n\text{ is odd,}\\
\left[
\begin{array}
[c]{lll}%
1 &  & \text{if }y=n-1,\\
2 &  & \text{otherwise}%
\end{array}
\right.  &  & \text{if }n\text{ is even.}%
\end{array}
\right.
\]
Consequently, the number of solus bitstrings of length $n$ with less than
$\ell$ $1$s and no runs of $k$ $0$s is%
\[
a_{n}=%
%TCIMACRO{\dsum \limits_{y=0}^{k-1}}%
%BeginExpansion
{\displaystyle\sum\limits_{y=0}^{k-1}}
%EndExpansion
\,%
%TCIMACRO{\dsum \limits_{x=0}^{\ell-1}}%
%BeginExpansion
{\displaystyle\sum\limits_{x=0}^{\ell-1}}
%EndExpansion
\,\widetilde{F}_{n}(n-x,y)
\]
and our prior results for $\ell=2,3,4,5$ and $k=7$ are easily verified. \ As
more examples, we have%
\begin{align*}
\{a_{n}\}_{n=1}^{6k-1}  &  =\left\{
2,3,5,8,13,21,33,52,83,132,209,327,502,752,1095,1543,2098,2749,\right. \\
&  \;\;\;\;\;\;3468,4210,4915,5517,5953,6173,6148,5876,5385,4727,3968,3178,\\
&  \left.  \;\;\;\;\;2422,1751,1196,767,458,252,126,56,21,6,1\right\}  ,
\end{align*}%
\begin{align*}
\{a_{n}\}_{n=1}^{7k-1}  &  =\left\{
2,3,5,8,13,21,33,52,83,132,210,334,530,836,1305,2005,3017,4428,\right. \\
&
\;\;\;\;\;\;6317,8739,11705,15163,18983,22957,26812,30236,32916,34582,35052,\\
&
\,\;\;\;\;\;34262,32277,29282,25556,21431,17242,13282,9772,6846,4550,2855,\\
&  \left.  \;\;\;\;\;1680,919,462,210,84,28,7,1\right\}
\end{align*}
for $\ell=6,7$ and $k=7$. \ If bitstrings are unconstrained, correlations for
$n=10,20,50$ turn out to be%
\[%
\begin{array}
[c]{ccccc}%
\rho(R_{10,0},S_{10})=-0.752444, &  & \rho(R_{20,0},S_{20})=-0.654958, &  &
\rho(R_{50,0},S_{50})=-0.530128.
\end{array}
\]
Dependency is more significant if bitstrings are solus:%
\[%
\begin{array}
[c]{ccccc}%
\rho(R_{10,0},S_{10})=-0.796825, &  & \rho(R_{20,0},S_{20})=-0.728540, &  &
\rho(R_{50,0},S_{50})=-0.616674.
\end{array}
\]
Table 2 exhibits values for larger $n=100,200,\ldots,1400$. \ While
acceleration of convergence is possible for each sequence, such techniques
merely suggest (without proof) that limits are nonzero as $n\rightarrow\infty
$. \ A rigorous method remains unknown.

\begin{center}%
\[%
\begin{tabular}
[c]{|c|c|c|}\hline
$n$ & $\rho$ for unconstrained case & $\rho$ for solus case\\\hline
$100$ & $-0.441772$ & $-0.525562$\\\hline
$200$ & $-0.361888$ & $-0.437637$\\\hline
$300$ & $-0.319761$ & $-0.389680$\\\hline
$400$ & $-0.292051$ & $-0.357617$\\\hline
$500$ & $-0.271797$ & $-0.333956$\\\hline
$600$ & $-0.256049$ & $-0.315434$\\\hline
$700$ & $-0.243295$ & $-0.300351$\\\hline
$800$ & $-0.232656$ & $-0.287715$\\\hline
$900$ & $-0.223581$ & $-0.276900$\\\hline
$1000$ & $-0.215704$ & $-0.267488$\\\hline
$1100$ & $-0.208773$ & $-0.259187$\\\hline
$1200$ & $-0.202606$ & $-0.251783$\\\hline
$1300$ & $-0.197066$ & $-0.245119$\\\hline
$1400$ & $-0.192050$ & $-0.239074$\\\hline
\end{tabular}
\ \ \ \ \ \
\]
Table 2:\ Correlation $\rho(R_{n,0},S_{n})$ as a function of $n$.\bigskip
\end{center}

Consider an unconstrained bitstring of length $n-1$. \ If we append the string
with a $1$, calling this $\sigma$, then there is a natural way \cite{HS-tcs4}
to associate $\sigma$ with an additive composition $\tau$ of $n$. \ For
example, if $n=10$,%
\[
\sigma=0110100111\longleftrightarrow\tau=\{2,1,2,3,1,1\}
\]
i.e., parts of $\tau$ correspond to \textquotedblleft waiting
times\textquotedblright\ for each $1$ in $\sigma$. \ The number of parts in
$\tau$ is equal to the bitsum $S_{n}$ of $\sigma$ and the maximum part in
$\tau$ is equal to the duration $R_{n,0}$ of the longest run of $0$s in
$\sigma$, plus one. \ Understanding correlation between these attributes,
given a uniformly distributed $\tau$, would seem vital. \ If instead we begin
with a solus $(n-1)$-bitstring starting and ending with $0$s, and append with
$1$ to construct $\sigma$, then the associated $\tau$ is a composition of $n$
with all parts $\geq2$. \ 

\section{Acknowledgements}

R, Mathematica and Maple have been useful throughout. I am indebted to a
friend, who wishes to remain anonymous, for giving encouragement and\ support
(in these dark days of the COVID-19 pandemic). \ I\ also recognize the editors
of the On-Line Encyclopedia of Integer Sequences for tireless and\ dedicated work.

\end{document}